\input amstex
\magnification\magstephalf
\documentstyle{amsppt}

\hsize 5.72 truein
\vsize 7.9 truein
\hoffset .39 truein
\voffset .26 truein
\mathsurround 1.67pt
\parindent 20pt
\normalbaselineskip 13.8truept
\normalbaselines
\binoppenalty 10000
\relpenalty 10000
\csname nologo\endcsname 


\font\bc=cmb10
\font\tenbsy=cmbsy10

\catcode`\@=11

\def\myitem#1.{\item"(#1)."\advance\leftskip10pt\ignorespaces}

\def\qedsymbol{{\mathsurround\z@$\square$}}
\redefine\qed{\relaxnext@\ifmmode\let\next\@qed\else
  {\unskip\nobreak\hfil\penalty50\hskip2em\null\nobreak\hfil
    \qedsymbol\parfillskip\z@\finalhyphendemerits0\par}\fi\next}
\def\@qed#1$${\belowdisplayskip\z@\belowdisplayshortskip\z@
  \postdisplaypenalty\@M\relax#1
  $$\par{\lineskip\z@\baselineskip\z@\vbox to\z@{\vss\noindent\qed}}}
\outer\redefine\beginsection#1#2\par{\par\penalty-250\bigskip\vskip\parskip
  \leftline{\tenbsy x\bf#1. #2}\nobreak\smallskip\noindent}
\outer\redefine\genbeginsect#1\par{\par\penalty-250\bigskip\vskip\parskip
  \leftline{\bf#1}\nobreak\smallskip\noindent}

\def\next{\let\@sptoken= }\def\next@{ }\expandafter\next\next@
\def\@futureletnext#1{\let\nextii@#1\futurelet\next\@flti}
\def\@flti{\ifx\next\@sptoken\let\next@\@fltii\else\let\next@\nextii@\fi\next@}
\expandafter\def\expandafter\@fltii\next@{\futurelet\next\@flti}

\let\zeroindent\z@
\let\savedef@\endproclaim\let\endproclaim\relax 
\define\chkproclaim@{\add@missing\endroster\add@missing\enddefinition
  \add@missing\endproclaim
  \envir@stack\endproclaim
  \edef\endit@{\leftskip\the\leftskip\rightskip\the\rightskip}}
\let\endproclaim\savedef@
\def\thing@{.\enspace\egroup\ignorespaces}
\def\thingi@(#1){ \rm(#1)\thing@}
\def\thingii@\cite#1{ \rm\@pcite{#1}\thing@}
\def\thingiii@{\ifx\next(\let\next\thingi@
  \else\ifx\next\cite\let\next\thingii@\else\let\next\thing@\fi\fi\next}
\def\thing#1#2#3{\chkproclaim@
  \ifvmode \medbreak \else \par\nobreak\smallskip \fi
  \noindent\advance\leftskip#1
  \hskip-#1#3\bgroup\bc#2\unskip\@futureletnext\thingiii@}
\let\savedef@\endproclaim\let\endproclaim\relax 
\def\endit{\endproclaim\endit@\let\endit@\undefined}
\let\endproclaim\savedef@
\def\defn#1{\thing\parindent{Definition #1}\rm}
\def\lemma#1{\thing\parindent{Lemma #1}\sl}
\def\prop#1{\thing\parindent{Proposition #1}\sl}
\def\thm#1{\thing\parindent{Theorem #1}\sl}
\def\cor#1{\thing\parindent{Corollary #1}\sl}

\def\remk#1{\thing\zeroindent{Remark #1}\rm}

\def\narrowthing#1{\chkproclaim@\medbreak\narrower\noindent
  \it\def\next{#1}\def\next@{}\ifx\next\next@\ignorespaces
  \else\bgroup\bc#1\unskip\let\next\narrowthing@\fi\next}
\def\narrowthing@{\@futureletnext\thingiii@}

\def\@cite#1,#2\end@{{\rm([\bf#1\rm],#2)}}
\def\cite#1{\in@,{#1}\ifin@\def\next{\@cite#1\end@}\else
  \relaxnext@{\rm[\bf#1\rm]}\fi\next}
\def\@pcite#1{\in@,{#1}\ifin@\def\next{\@cite#1\end@}\else
  \relaxnext@{\rm([\bf#1\rm])}\fi\next}

\advance\minaw@ 1.2\ex@
\atdef@[#1]{\ampersand@\let\@hook0\let\@twohead0\brack@i#1,\z@,}
\def\brack@{\z@}
\let\@@hook\brack@
\let\@@twohead\brack@
\def\brack@i#1,{\def\next{#1}\ifx\next\brack@
  \let\next\brack@ii
  \else \expandafter\ifx\csname @@#1\endcsname\brack@
    \expandafter\let\csname @#1\endcsname1\let\next\brack@i
    \else \Err@{Unrecognized option in @[}%
  \fi\fi\next}
\def\brack@ii{\futurelet\next\brack@iii}
\def\brack@iii{\ifx\next>\let\next\brack@gtr
  \else\ifx\next<\let\next\brack@less
    \else\relaxnext@\Err@{Only < or > may be used here}
  \fi\fi\next}
\def\brack@gtr>#1>#2>{\setboxz@h{$\m@th\ssize\;{#1}\;\;$}%
 \setbox@ne\hbox{$\m@th\ssize\;{#2}\;\;$}\setbox\tw@\hbox{$\m@th#2$}%
 \ifCD@\global\bigaw@\minCDaw@\else\global\bigaw@\minaw@\fi
 \ifdim\wdz@>\bigaw@\global\bigaw@\wdz@\fi
 \ifdim\wd@ne>\bigaw@\global\bigaw@\wd@ne\fi
 \ifCD@\enskip\fi
 \mathrel{\mathop{\hbox to\bigaw@{$\ifx\@hook1\lhook\mathrel{\mkern-9mu}\fi
  \setboxz@h{$\displaystyle-\m@th$}\ht\z@\z@
  \displaystyle\m@th\copy\z@\mkern-6mu\cleaders
  \hbox{$\displaystyle\mkern-2mu\box\z@\mkern-2mu$}\hfill
  \mkern-6mu\mathord\ifx\@twohead1\twoheadrightarrow\else\rightarrow\fi$}}%
 \ifdim\wd\tw@>\z@\limits^{#1}_{#2}\else\limits^{#1}\fi}%
 \ifCD@\enskip\fi\ampersand@}
\def\brack@less<#1<#2<{\setboxz@h{$\m@th\ssize\;\;{#1}\;$}%
 \setbox@ne\hbox{$\m@th\ssize\;\;{#2}\;$}\setbox\tw@\hbox{$\m@th#2$}%
 \ifCD@\global\bigaw@\minCDaw@\else\global\bigaw@\minaw@\fi
 \ifdim\wdz@>\bigaw@\global\bigaw@\wdz@\fi
 \ifdim\wd@ne>\bigaw@\global\bigaw@\wd@ne\fi
 \ifCD@\enskip\fi
 \mathrel{\mathop{\hbox to\bigaw@{$%
  \setboxz@h{$\displaystyle-\m@th$}\ht\z@\z@
  \displaystyle\m@th\mathord\ifx\@twohead1\twoheadleftarrow\else\leftarrow\fi
  \mkern-6mu\cleaders
  \hbox{$\displaystyle\mkern-2mu\copy\z@\mkern-2mu$}\hfill
  \mkern-6mu\box\z@\ifx\@hook1\mkern-9mu\rhook\fi$}}%
 \ifdim\wd\tw@>\z@\limits^{#1}_{#2}\else\limits^{#1}\fi}%
 \ifCD@\enskip\fi\ampersand@}


\define\ie{{\it i.e\.}}
\define\today{\number\day\ \ifcase\month\or
  January\or February\or March\or April\or May\or June\or
  July\or August\or September\or October\or November\or December\fi
  \ \number\year}
\def\pr@m@s{\ifx'\next\let\nxt\pr@@@s \else\ifx^\next\let\nxt\pr@@@t
  \else\let\nxt\egroup\fi\fi \nxt}

\define\widebar#1{\mathchoice
  {\setbox0\hbox{\mathsurround\z@$\displaystyle{#1}$}\dimen@.1\wd\z@
    \ifdim\wd\z@<.4em\relax \dimen@ -.16em\advance\dimen@.5\wd\z@ \fi
    \ifdim\wd\z@>2.5em\relax \dimen@.25em\relax \fi
    \kern\dimen@ \overline{\kern-\dimen@ \box0\kern-\dimen@}\kern\dimen@}%
  {\setbox0\hbox{\mathsurround\z@$\textstyle{#1}$}\dimen@.1\wd\z@
    \ifdim\wd\z@<.4em\relax \dimen@ -.16em\advance\dimen@.5\wd\z@ \fi
    \ifdim\wd\z@>2.5em\relax \dimen@.25em\relax \fi
    \kern\dimen@ \overline{\kern-\dimen@ \box0\kern-\dimen@}\kern\dimen@}%
  {\setbox0\hbox{\mathsurround\z@$\scriptstyle{#1}$}\dimen@.1\wd\z@
    \ifdim\wd\z@<.28em\relax \dimen@ -.112em\advance\dimen@.5\wd\z@ \fi
    \ifdim\wd\z@>1.75em\relax \dimen@.175em\relax \fi
    \kern\dimen@ \overline{\kern-\dimen@ \box0\kern-\dimen@}\kern\dimen@}%
  {\setbox0\hbox{\mathsurround\z@$\scriptscriptstyle{#1}$}\dimen@.1\wd\z@
    \ifdim\wd\z@<.2em\relax \dimen@ -.08em\advance\dimen@.5\wd\z@ \fi
    \ifdim\wd\z@>1.25em\relax \dimen@.125em\relax \fi
    \kern\dimen@ \overline{\kern-\dimen@ \box0\kern-\dimen@}\kern\dimen@}%
  }

\catcode`\@\active

\let\PVstyle=d 

\input xy
\xyoption{matrix}
\xyoption{arrow}
\xyoption{cmtip}
\xyoption{dvips}

\font\tenscr=rsfs10 
\font\sevenscr=rsfs7 
\font\fivescr=rsfs5 
\skewchar\tenscr='177 \skewchar\sevenscr='177 \skewchar\fivescr='177
\newfam\scrfam \textfont\scrfam=\tenscr \scriptfont\scrfam=\sevenscr
\scriptscriptfont\scrfam=\fivescr
\define\scr#1{{\fam\scrfam#1}}
\let\Cal\scr

\let\0\relax 
\def\restrictedto#1{\big|_{#1}}
\define\tildeto{\mathrel{\mathsurround0pt\rlap{\raise.8ex\hbox{\kern.06em%
  $\scriptstyle\sim$}}\hbox{$\to$}}}

\define\pr{\operatorname{pr}}
\define\Proj{\operatorname{Proj}}

\define\Spec{\operatorname{Spec}}

\topmatter
\title Nagata's embedding theorem\endtitle
\author Paul Vojta\endauthor
\affil University of California, Berkeley\endaffil
\address Department of Mathematics, University of California,
  Berkeley, CA \ 94720\endaddress
\date 13 June 2007 \enddate
\thanks Supported by NSF grants DMS95-32018 and DMS-0500512, and
the Institute for Advanced Study.\endthanks
\subjclassyear{2000}

\abstract
In 1962--63, M. Nagata showed that an abstract variety could be embedded
into a complete variety.  Later, P. Deligne translated Nagata's proof
into the language of schemes, but did not publish his notes.  This paper,
which is to appear as an appendix in a forthcoming book, gives an elaboration
of Deligne's notes.  It also contains some complementary results on extending
divisors and vector sheaves to suitable completions.\endabstract
\endtopmatter

\document

The goal of this note is to prove that, if $X$ is a scheme,
separated and of finite type over a noetherian scheme $S$, then there
exists a proper $S$\snug-scheme $\widebar X$ and an open immersion
$X\hookrightarrow\widebar X$ over $S$ with schematically dense image
(Theorem \04.1).  In addition, given a Cartier divisor $D$ or a
vector sheaf $\Cal E$ on $X$, the completion $\widebar X$ can be chosen
so that $D$ or $\Cal E$ extends to a Cartier divisor or vector sheaf
(respectively) on $\widebar X$.

The first assertion was proved by Nagata; see \cite{N~1}, \cite{N~2}.
Nagata's proof is phrased in terms of Zariski's language of algebraic geometry,
though, which makes it difficult for many to read.  Because of this,
P. Deligne wrote some notes \cite{D}, which translate Nagata's work
into the language of schemes.  These notes, however, are unpublished.

This paper, which is based closely on Deligne's notes, was written
to give a mostly self-contained exposition of the proof.
After writing these notes, I encountered another (much more
thorough) rendition of Deligne's notes by B. Conrad \cite{C}.
Conrad also notes the existence of another proof of Nagata's theorem
by L\"utkebohmert \cite{L}, which uses schemes but uses different methods.

Sections \01--\04 of this note give Nagata's proof, following \cite{D}.
Section \05 adds some complementary results on constructing the completion
so that certain sheaves or divisors extend to the completion.

{\smallskip\narrower\it  Throughout this note all schemes are assumed
to be noetherian and all ideal sheaves to be coherent.
Projective morphisms are as defined in \cite{EGA}.
Schemes are not assumed to be
separated unless it is explicitly mentioned.\par\smallskip}

Given an open subscheme $U$ of a scheme $X$, there are two notions of $U$
being dense in $X$, and in this note the difference between them will
be important.  The weaker notion is that the underlying topological space
of $U$ is dense in the topological space of $X$.  The wording ``$U$ is
a dense open subset of $X$'' or ``$U$ is Zariski-dense in $X$'' will describe
this weaker notion.  A stronger notion is that the scheme-theoretic image
of $U$ under the open immersion $U\hookrightarrow X$ is all of $X$;
equivalently, $U$ contains all associated points of $X$\footnote{This is shown
by reducing to the affine case and using the primary decomposition.}.
The wording ``$U$ is an open dense subscheme of $X$'' or
``$U$ is schematically dense in $X$'' will refer to this notion.
This stronger notion of denseness is needed for Lemma \02.2e (for example),
which is used in the proof of Theorem \02.4.

\beginsection{\01}{Preliminary results on blowings-up}

We start with some lemmas on blowings-up.  For the definition of the
blowing-up of a scheme $X$ along an ideal sheaf in $\Cal O_X$,
see \cite{H, II~\S\,7}.

\lemma{\01.1}  Let $X$ be a scheme and let $\frak a_1,\dots,\frak a_n$
be ideal sheaves in $\Cal O_X$.  Let $\widetilde X$ be the blowing-up
of $X$ along $\frak a_1+\dots+\frak a_n$.  Then the strict transforms
of $V(\frak a_1),\dots,V(\frak a_n)$ in $\widetilde X$ have empty
intersection.
\endit

\demo{Proof}  The case $n=2$ is \cite{H, II Ex.~7.12}; the general case
is analogous.\qed
\enddemo

\lemma{\01.2}  Let $U$ be an open subscheme of a scheme $X$
and let $\frak a_0$ be an ideal sheaf on $U$.  Then there exists a (coherent)
ideal sheaf $\frak a$ on $X$ such that $\frak a\restrictedto U=\frak a_0$
and such that $V(\frak a)=\widebar{V(\frak a_0)}$.
\endit

\demo{Proof}  This follows from \cite{EGA, I~9.5.10}
and \cite{H, II Proposition~5.9}.\qed
\enddemo

\lemma{\01.3}  Let $X$ be a scheme and let $\frak a$ be an ideal
sheaf in $\Cal O_X$.  Let $\pi\:\widetilde X\to X$ be the blowing-up of $X$
along $\frak a$.  Then
$$\pi^{-1}(V(\frak a)) = V(\pi^{-1}\frak a\cdot\Cal O_{\widetilde X})$$
as closed subsets of $\widetilde X$.
\endit

\demo{Proof}  See \cite{EGA, II~8.1.8}.\qed
\enddemo

\lemma{\01.4}  Let $X$ be a scheme, and let $\frak a$ and $\frak b$ be
ideal sheaves in $\Cal O_X$.  Let $\pi_1\:X_1\to X$ and $\pi_2\:X_2\to X$ be
the blowings-up of $X$ along $\frak a$ and $\frak a\frak b$, respectively.
Then $X_2$ dominates $X_1$ (\ie, $\pi_2$ factors uniquely through $\pi_1$).
\endit

\demo{Proof}  The question is local on $X$, so we may assume that $X$
is affine, say $X=\Spec A$.  Then we may regard $\frak a$ and $\frak b$
as ideals in $A$.  Let $a_0,\dots,a_n$ and $b_0,\dots,b_m$ be systems
of generators for $\frak a$ and $\frak b$, respectively.
Recall that $X_1=\Proj S$, where $S$ is the graded ring
$\bigoplus_{d\ge0}\frak a^d$ (where $\frak a^0=A$).  This is covered
by open affines $\Spec S_{(a_i)}$ for $i=0,\dots,n$.  These are glued
by identifying $D(a_j/a_i)$ in $\Spec S_{(a_i)}$
with $D(a_i/a_j)$ in $\Spec S_{(a_j)}$ via the canonical isomorphism
$$\bigl(S_{(a_i)}\bigr)_{a_j/a_i} \cong \bigl(S_{(a_j)}\bigr)_{a_i/a_j}.$$
The other blowing-up $X_2$ is similarly covered by open affines
$\Spec T_{(a_ib_j)}$, where $T=\bigoplus_{d\ge0}\frak a^d\frak b^d$.
The ring $S_{(a_i)}$ is generated over $A$ by $a_0/a_i,\dots,a_n/a_i$;
thus the association $a_\ell/a_i\mapsto a_\ell b_j/a_ib_j$ defines a ring
homomorphism $\phi_{ij}\:S_{(a_i)}\to T_{(a_ib_j)}$ for all $i$ and $j$.
These glue to give a morphism $X_2\to X_1$.  This morphism commutes with
the maps to $X$ because all $\phi_{ij}$ commute with the maps $A\to S_{(a_i)}$
and $T_{(a_ib_j)}$.\qed
\enddemo

\prop{\01.5} \cite{R-G, Pr\`emiere partie, Lemme 5.1.4}  Let $\pi_1\:X_1\to X$
be the blowing-up of a scheme $X$ along an ideal sheaf $\frak a$ in $\Cal O_X$,
and let $\pi_2\:X_2\to X_1$ be the blowing-up of $X_1$ along an ideal
sheaf $\frak b$ in $\Cal O_{X_1}$.  Then there exists an ideal sheaf $\frak c$
in $\Cal O_X$ such that the blowing-up $\pi\:\widetilde X\to X$ of $X$
along $\frak c$ is isomorphic to $\pi_2\circ\pi_1\:X_2\to X$.
Moreover, $\frak c$ can be chosen such
that $V(\frak c)=V(\frak a)\cup \pi_1(V(\frak b))$.
\endit

\demo{Proof}  {\bc Case I.}  $X$ is affine.  Say $X=\Spec A$.
Then $\frak a$ may be regarded as an ideal in $A$.

By \cite{EGA, II~4.4.3 and II~4.5.10}, the line sheaf $\Cal O(1)$
on $X_1$ (defined via the blowing-up $\pi_1$) is ample.  Hence,
for all sufficiently large integers $n$, $\frak b(n):=\frak b\otimes\Cal O(n)$
is generated by global sections.  Since $\frak b\subseteq\Cal O_{X_1}$,
these global sections lie in $\Gamma(X_1,\Cal O(n))$.

Let $x_0,\dots,x_r$ be a system of generators for $\frak a$.
These determine a graded homomorphism
$\phi\:A[X_0,\dots,X_r]\to\bigoplus_{d\ge0}\frak a^d$, which in turn defines
a closed immersion $f\:X_1\to\Bbb P^r_A$.  Let $I=\ker\phi$; we then have an
exact sequence of sheaves on $\Bbb P^r_A$:
$$0 \to \widetilde I \to \Cal O_{P^r_A} \to f_{*}\Cal O_{X_1} \to 0.$$
Tensor with $\Cal O(n)$ and take the long exact sequence in cohomology.
By Serre's theorem \cite{H, III~5.2}, $H^1(\Bbb P^r_A,\widetilde I(n))=0$
for all $n$ sufficiently large; for these $n$, $\Gamma(\Bbb P^r_A,\Cal O(n))$
maps surjectively onto $\Gamma(X_1,\Cal O(n))$.  Thus, for $n$ sufficiently
large, $\Gamma(X_1,\Cal O(n))=\frak a^n$.

Now fix an $n$ such that $\frak b(n)$ is generated by global sections
and such that $\Gamma(X_1,\Cal O(n))=\frak a^n$.  Let $b_0,\dots,b_m$
be a set of global sections that generate $\frak b(n)$; by what was shown
earlier, these may be regarded as elements of $\frak a^n$.

Let $\frak c$ be the ideal sheaf on $X$ corresponding to the ideal
$$\frak a^n(b_0,\dots,b_m)\subseteq A,$$
and let $\widetilde X$ be the blowing-up of $X$ along $\frak c$.
We claim that $\widetilde X$ is isomorphic to $X_2$.
Indeed, let $a_0,\dots,a_s$ be a system of generators for $\frak a^n$,
and let $S=\bigoplus_{d\ge0}\frak a^{nd}$.  Then $X_1=\Proj S$;
it is covered by open affines $\Spec S_{(a_i)}$, $i=0,\dots,s$.
For each such $i$, the ideal $\frak b_i\subseteq S_{(a_i)}$ corresponding
to the ideal sheaf $\frak b$ equals
$$\frak b_i = \left(\frac{b_0}{a_i},\dots,\frac{b_m}{a_i}\right).$$
Let $T_i=\bigoplus_{d\ge0}\frak b_i^d$;
then $\pi_2^{-1}(\Spec S_{(a_i)})=\Proj T_i$; it is covered by
open affines $\Spec (T_i)_{(b_j/a_i)}$, $j=0,\dots,m$.
Also let $S'=\bigoplus_{d\ge0}\frak c^d$; then $\widetilde X=\Proj S'$
is covered by open affines $\Spec S'_{(a_ib_j)}$ for $i=0,\dots,s$ and
$j=0,\dots,m$.  The result then follows, in this case, since
$$(T_i)_{(b_j/a_i)}\cong S'_{(a_ib_j)}.$$
(Verification of this isomorphism is left as an exercise for the reader.)

By construction it is clear that $V(\frak c)=V(\frak a)\cup \pi_1(V(\frak b))$.

{\bc Case II.}  General case.  Let $U_1,\dots,U_\ell$ be a cover of $X$
by open affines.  Let $\frak c_1,\dots,\frak c_\ell$ be the ideal sheaves
constructed as in Case I, using the same value of $n$.  Since the construction
in Case I commutes with localizing to the local ring at a point of $\Spec A$,
these ideal sheaves glue to give an ideal sheaf $\frak c$ on $X$
with the desired properties.\qed
\enddemo

\comment
\lemma{\01.6}  Let $i\:X'\to X$ be a closed immersion of schemes,
let $\frak a'$ be an ideal sheaf on $X'$, and let $\frak a=i_{*}\frak a'$.
Let $\widetilde X'$ (resp\. $\widetilde X$) be the blowings-up of $X'$
(resp\. $X$) along $\frak a'$ (resp\. $\frak a$), respectively.
Then there is a unique morphism $\tilde\imath\:\widetilde X'\to\widetilde X$
making the diagram
$$\CD \widetilde X' @>\tilde\imath>> \widetilde X \\
  @VVV @VVV \\
  X' @>i>> X \endCD$$
commute; moreover, $\tilde\imath$ is a closed immersion.
\endit

\demo{Proof}  Indeed, $\frak a'=i^{-1}\frak a\cdot\Cal O_{X'}$,
and the result then follows from \cite{H, II~7.15}.\qed
\enddemo
\endcomment

\lemma{\01.6}  Let $X$ be a scheme, let $U$ be an open dense subscheme,
and let $\frak a$ be an ideal sheaf on $X$ with $V(\frak a)$ disjoint from $U$.
Let $\pi\:\widetilde X\to X$ denote the blowing-up of $X$ along $\frak a$.
Then $\pi^{-1}(U)$ is schematically dense in $\widetilde X$.
\endit

\demo{Proof}  This holds because, by \cite{H, II~7.13a}, the inverse image
ideal sheaf of $\frak a$ in $\widetilde X$ is a line sheaf.  Therefore
the exceptional divisor is Cartier, so locally it is generated by a function
that is not a zero divisor.  Thus the exceptional divisor does not contain
any associated points.\qed
\enddemo

\lemma{\01.7}  Let $X$ be an open subscheme of a scheme $\widebar X$,
let $p\:\widebar Z\to \widebar X$ be a morphism of finite type,
let $Z=p^{-1}(X)$, and let $F$ and $G$ be closed subsets of $X$ and $Z$,
respectively, with $p^{-1}(F)\cap G=\emptyset$.
$$\CD G &\;\;\subseteq\;\;& Z @[hook]>>> \widebar Z \\
  && @VVV @VVpV \\
  F &\;\;\subseteq\;\;& X @[hook]>>> \widebar X \endCD$$
Then there exists an ideal sheaf $\frak a$ on $\widebar X$, with $V(\frak a)$
disjoint from $X$, with the following property.  Let $\widetilde X$ denote
the blowing-up of $\widebar X$ along $\frak a$, let
$\tilde p\:\widetilde Z\to\widetilde X$ be the morphism obtained from $p$
by base change, let $\widetilde F$ denote the closure of $F$ in $\widetilde X$
(where $X$ is regarded as an open subscheme of $\widetilde X$), and likewise
let $\widetilde G$ be the closure of $G$ in $\widetilde Z$.  Then
$$\tilde p^{-1}(\widetilde F) \cap \widetilde G = \emptyset.$$
\endit

\demo{Proof}  We show that this holds for any ideal sheaf $\frak a$ satisfying:
\roster
\myitem i.  $V(\frak a)$ is disjoint from $X$;
\myitem ii.  $\frak a$ contains the ideal corresponding to $\widebar F$; and
\myitem iii.  $\frak a$ is contained in the ideal corresponding to the
scheme-theoretic image $p(p^{-1}(\widebar F)\cap\widebar G)$.
\endroster
Such ideals exist:  take $\frak a$ equal to the ideal in (iii), for example.

The above conditions are preserved upon passing to open subsets of $\widebar X$
and $\widebar Z$, so we may assume that $\widebar X$ and $\widebar Z$ are
affine, say $\widebar X=\Spec A$ and $\widebar Z=\Spec B$.  Let $\phi\:A\to B$
be the ring homomorphism corresponding to $p$.  By abuse of notation,
let $\frak a$ refer now to an ideal in the ring $A$.  Also let $\frak f$
and $\frak g$ denote the ideals in $A$ and $B$ corresponding to the
closed subsets $\widebar F$ and $\widebar G$, respectively.

In the language of rings, the counterparts to conditions (ii) and (iii)
above are
$$\frak a\supseteq\frak f\tag\01.7.1$$
and
$$\frak a\subseteq \phi^{-1}(\phi(\frak f)B+\frak g).\tag\01.7.2$$
Let $f_0,\dots,f_r$ be a system of generators for $\frak f$.  As $x$ varies
over a system of generators for $\frak a$, the sets $\Spec S_{(x)}$
give an open affine cover of $\widetilde X$,
where $S=\bigoplus_{i\ge0}\frak a^i$.  By (\01.7.2)
we have $\phi(x)=\sum_{i=0}^r b_i\phi(f_i)+g$ for some $b_0,\dots,b_r\in B$
and $g\in\frak g$.  We have $g\in\phi(\frak a)B$ since
$g=\phi(x)-\sum b_i\phi(f_i)$, $\phi(x)\in\phi(\frak a)B$ (by choice of $x$),
and $\phi(f_i)\in\phi(\frak a)$ by (\01.7.1).  Therefore we may write
$$1 = \frac xx = \frac{\sum b_i\phi(f_i)}{x} + \frac gx$$
in $B\otimes_A S_{(x)}$.  But the first term lies in the ideal of
$\tilde p^{-1}(\widetilde F)$, and the second in the ideal of $\widetilde G$;
hence these two closed sets are disjoint, as was to be shown.\qed
\enddemo

Note that the set $p(G)$ is not necessarily closed, so this cannot be proved
by applying Lemma \01.1 to $\widebar F$ and $p(\widebar G)$.

\beginsection{\02}{Quasi-dominations and extensions thereof}

This section defines quasi-dominations and shows that they can be extended
to a larger domain scheme after blowing up.

\defn{\02.1}  Let $S$ be a noetherian scheme, and let $X$ and $Y$ be
separated $S$\snug-schemes of finite type.
\roster
\myitem a.  A {\bc quasi-domination} $f\:X\dashrightarrow Y$ is a pair
consisting of an open dense subscheme $U\subseteq X$ and a morphism
$f\:U\to Y$ whose graph is closed in $X\times_S Y$.
\myitem b.  Such a quasi-domination is {\bc proper} if its underlying
morphism $f\:U\to Y$ is proper.
\myitem c.  If $f\:X\dashrightarrow Y$ is a quasi-domination (resp\. a
proper quasi-domination), then we say that $X$ {\bc quasi-dominates} $Y$
(resp\. {\bc properly quasi-dominates} $Y$) via $f$.  Mention of $f$ may
be omitted if it is clear from the context.
\endroster
\endit

\lemma{\02.2} Let $S$, $X$, and $Y$ be as above.
\roster
\myitem a.  Let $\psi\:U\to Y$ be a quasi-domination $X\dashrightarrow Y$,
and let $V$ be an open subset of $Y$ such that $\psi^{-1}(V)$ is schematically
dense in $X$.  Then $\psi\restrictedto{\psi^{-1}(V)}$ determines a
quasi-domination $X\dashrightarrow V$.
\myitem b.  If $X$ quasi-dominates $Y$ and if $U$ is an open subset of $X$,
then $U$ quasi-dominates $Y$.
\myitem c.  Let $f\:X\to S$ denote the structural morphism of $X$,
let $W$ be an open subset of $S$, and suppose that $Y$ is a scheme over $W$.
If $f^{-1}(W)$ is schematically dense in $X$, then any quasi-domination
$f^{-1}(W)\dashrightarrow Y$ is also a quasi-domination $X\dashrightarrow Y$.
\myitem d.  Let $\psi\:X\dashrightarrow Y$ be a quasi-domination,
let $\frak a$ be an ideal sheaf on $X$ such that $X\setminus V(\frak a)$
is schematically dense, and let $\pi\:\widetilde X\to X$ be the blowing-up
of $X$ along $\frak a$.  Then $\psi\circ\pi$ is a quasi-domination
$\widetilde X\dashrightarrow Y$.
\myitem e.  If $U\subseteq X$ is schematically dense and $\psi_1\:X\to Y$
and $\psi_2\:X\to Y$ are morphisms such that
$\psi_1\restrictedto U=\psi_2\restrictedto U$, then $\psi_1=\psi_2$.
\myitem f.  If $U$ is an open dense subscheme of $X$ and $f\:U\to Y$ is
proper, then $f$ is a proper quasi-domination $X\dashrightarrow Y$.
\endroster
\endit

\demo{Proof}  Part (a) follows from the fact that the intersection of the
graph of $\psi$ with $X\times_S V$ will be a closed subset of $X\times_S V$,
and from the fact that the domain will still be schematically dense.
Part (b) is similar.
Part (c) follows from the fact that $f^{-1}(W)\times_W Y=X\times_S Y$.
Part (d) follows from Lemma \01.6 and the fact that the graph of $\psi\circ\pi$
is the pull-back via $(\pi\times_S 1)$ of the graph of $\psi$.
Part (e) follows from the fact that the inverse image of the diagonal
by the morphism $(f,g)\:X\to Y\times_S Y$ is a closed subscheme of $X$
containing $U$, hence is all of $X$.
Finally, (f) follows from the fact that the graph of the inclusion
$U\hookrightarrow X$ is a closed subset of $U\times_S X$, and that the map
$U\times_S X\to Y\times_S X$ is proper.\qed
\enddemo

\lemma{\02.3}  Let $X$ be a separated $S$\snug-scheme of finite type
and let $U$ be an open dense subscheme of $X$.
Let $Y$ be a closed subscheme of $\Bbb A^r_S$ for some $r\in\Bbb N$,
and let $\phi\:Y\hookrightarrow P$ be its scheme-theoretic closure
in $\Bbb P^r_S$.  Let $\psi\:U\to Y$ be a morphism over $S$.
Then there exists an ideal sheaf $\frak a$ on $X$,
with $V(\frak a)$ disjoint from $U$, such that if $\pi\:\widetilde X\to X$
is the blowing-up of $X$ along $\frak a$, then $\psi$ extends to a
morphism $\widetilde\psi\:\widetilde X\to P$:
$$\CD U @[hook]> \pi^{-1}\restrictedto U >> \widetilde X \\
  @VV \psi V @VV \widetilde\psi V \\
  Y @[hook]> \phi >> P \endCD$$
\endit

\demo{Proof} Let $x_1,\dots,x_r$ be the standard coordinate functions
on $\Bbb A^r_S$.

{\bc Case I.}  The scheme $X$ is affine, say $X=\Spec A$,
and $U=D(u)$ for some $u\in A$ which is not a zero divisor.
Then, since $\psi$ is a morphism, $\psi^{*}x_1,\dots,\psi^{*}x_r$ lie in
the localized ring $A_u$; hence there exists a large integer $n$ such that
$$u^n\psi^{*}x_1,\dots,u^n\psi^{*}x_r\in A.$$
Let
$$\frak a = \bigl(u^n,u^n\psi^{*}x_1,\dots,u^n\psi^{*}x_r\bigr).$$
Then the blowing-up $\widetilde X$ of $X$ along $\frak a$ admits an
extension of $\psi$ to a morphism $\widetilde X\to P$.

{\bc Case II.}  The complement $X\setminus U$ is the support of an effective
Cartier divisor $D$.  Let $\{W_1,\dots,W_m\}$ be an open affine cover of
$X$ such that $D\restrictedto{W_i}$ is principal for each $i$.  By Case I
there exists an ideal sheaf $\frak a_i$ on $W_i$ for each $i$ which works.
These extend to ideal sheaves $\tilde{\frak a}_i$ on $X$ for each $i$
by Lemma \01.2.  Then by Lemma \01.4 it suffices to blow up the product of
the $\tilde{\frak a}_i$.  This gives morphisms $W_i\to P$ which extend
$\psi\restrictedto{W_i\cap U}\:W_i\cap U\to P$.  Since $U$ is schematically
dense in $X$, it is schematically dense in $\widetilde X$; hence these
extensions are unique and therefore glue to give a morphism
$\widetilde\psi\:\widetilde X\to P$.

{\bc Case III.}  The lemma holds in general.  Indeed, one can reduce to
Case II by blowing up the ideal sheaf defining the closed subset
$X\setminus U$.\qed
\enddemo

\thm{\02.4} (Deligne \cite{D}; cf\. \cite{N~1, Theorem~3.2})
Let $X$ and $Y$ be separated $S$\snug-schemes of finite type, let $U$ be
an open dense subscheme of $X$, and let $\psi\:U\to Y$ be a morphism
over $S$.  Then there exists an ideal sheaf $\frak a$ on $X$,
with $V(\frak a)$ disjoint from $U$, such that if $\widetilde X$
is the blowing-up of $X$ along $\frak a$, then $\psi$ extends to a
quasi-domination $\widetilde\psi\:\widetilde X\dashrightarrow Y$.
\endit

\demo{Proof}  Since $Y$ is of finite type over $S$, there is an open
cover $\{Y_1,\dots,Y_n\}$ of $Y$ by open affines such that the affine
ring of each $Y_i$ is an algebra of finite type over the affine ring of some
open $W_i\subseteq S$.
For each $i$ let $U_i=\psi^{-1}(Y_i)$, let $F_i=U\setminus U_i$,
and let $\Gamma_i^0\subseteq U_i\times_S Y_i$ denote the graph of
$\psi\restrictedto{U_i}$.  By Lemma \01.7 we may further blow up $X$,
not touching $U$, such that the closure $\widebar\Gamma_i^0$ of $\Gamma_i^0$
in $X\times_S Y_i$ does not touch $\widebar F_i\times_S Y_i$.
Let $X_i=X\setminus\widebar F_i$.

For each $i$, apply Lemma \02.3 to $U_i\subseteq X_i$ and to
$\psi\restrictedto{U_i}\:U_i\to Y_i$.  This gives a projective completion
$\phi_i\:Y_i\hookrightarrow P_i$ over $W_i$ and an ideal sheaf $\frak a_i$
on $X_i$, disjoint from $U_i$, such that after blowing up, we have a morphism
$\psi_i\:X_i\to P_i$ extending $\phi_i\circ\psi\restrictedto{U_i}$.
In particular, by Lemma \02.2a and Lemma \02.2c, $\psi_i$ determines
a quasi-domination $X_i\dashrightarrow Y_i$ for each $i$.
By Lemmas \01.2, \01.4, and \02.2d, there is a blowing-up
of $X$ along an ideal sheaf $\frak a$ with $V(\frak a)\cap U=\emptyset$
such that $\psi\restrictedto{U_i}$ extends to a quasi-domination
$X_i\dashrightarrow Y_i$ simultaneously for all $i$.
For each $i$ let $\Gamma_i$ denote the graph of this quasi-domination;
it is a closed subset of $X_i\times_S Y_i$.
Let $\Gamma\subseteq X\times_S Y$ be the union of the $\Gamma_i$.

The projection $X\times_S Y\to X$ induces isomorphisms $\Gamma_i\tildeto V_i$
for some open subsets $V_i\subseteq X_i$.  For all $i$ and $j$ the
subscheme $V_i\cap V_j\cap U$ is dense in $V_i\cap V_j$;
therefore $\Gamma_i$ and $\Gamma_j$ coincide over $V_i\cap V_j$ by
Lemma \02.2e.  Thus the schemes $\Gamma_i$ glue
to give $\Gamma$ the structure of a scheme such that the projection
$X\times_S Y\to X$ induces an isomorphism $\Gamma\tildeto\bigcup V_i$.
Hence $\Gamma$ is the graph of a morphism $\widetilde\psi\:\bigcup V_i\to Y$.

It remains to show that $\Gamma$ is a closed subset.
Suppose that $\Gamma$ is not closed; then there exists a point
$\xi\in\widebar\Gamma$ such that $\xi\notin\Gamma$.  There exists an
irreducible component $X'$ of $X$ and an index $i$ such that
$\xi\in\widebar{\Gamma_i\cap(X'\times_S Y)}$.  Let $\eta$ denote the
generic point of $X'$; then $\xi\in\widebar{\{(\eta,\psi(\eta))\}}$.
Pick an index $j$ such that $\xi\in X\times_S Y_j$.
We have $Y_j\cap\widebar{\{\psi(\eta)\}}\ne\emptyset$, so $\psi(\eta)\in Y_j$,
and therefore $\eta\in X_j$.
It follows that $\widebar{\{(\eta,\psi(\eta))\}}\subseteq\widebar\Gamma_j$
and therefore that $\xi\in\widebar\Gamma_j$.  But $\Gamma_j$ is a closed
subset of $X_j\times_S Y_j$ and $\xi\in X\times_S Y_j$;
therefore $\xi\in\widebar F_j\times_S Y_j$.  Thus the closure of $\Gamma_j$
in $X\times_S Y_j$ meets $\widebar F_j\times_S Y_j$, a contradiction.

Thus $\Gamma$ determines a quasi-domination $X\dashrightarrow Y$,
as was to be shown.\qed
\enddemo

\cor{\02.5} (Chow; \cite{D, Cor.~1.4})  Let $X$ be a separated $S$\snug-scheme
of finite type.  Let $U$ be an open dense subset of $X$, quasi-projective
over $S$.  Then there exists a diagram
$$\xymatrix@R=4mm{ & X' \ar@{^{(}->}[]+<4mm,0mm>;[r] \ar[dd]^q
    & \widebar X \\ 
  U \ar@{^{(}->}[]+<4mm,2mm>;[ur] \ar@{^{(}->}[dr] \\ 
  & X }\tag\02.5.1$$
in which $q$ is a proper morphism, isomorphic over $U$;
$X\hookrightarrow\widebar X$ is an open immersion; and $\widebar X$ is
projective over $S$.  We may also assume that the image of $U$ is
schematically dense in $\widebar X$ (and therefore in $X'$).
\endit

\demo{Proof}  Let $W$ be a projective $S$\snug-scheme containing $U$ as
an open dense subscheme.  By Theorem \02.4 applied to the open immersion
$U\hookrightarrow X$, we may assume that $W$ quasi-dominates $X$.
We may then take $\widebar X=W$.  Indeed, let $X'$ denote the domain
of the quasi-domination $W\dashrightarrow X$.  The graph of this
quasi-domination is closed in $W\times_S X$, and $W\times_S X$ is proper
over $X$, so this graph, and hence $X'$, is proper over $X$.

Projectivity of $\widebar X$ follows from \cite{EGA, II 5.5.5, ii}.\qed
\enddemo

The special case in which $X$ is proper over $S$ is important enough to state
separately.  See also \cite{R-G, Pr\`emiere partie, Cor.~5.7.14}.

\cor{\02.6}  Let $X$ be a proper scheme over $S$, and let $U$ be an open
dense subset of $X$.  If $U$ is quasi-projective over $S$, then
there is a proper morphism $q\:X'\to X$, isomorphic over $U$,
such that $X'$ is projective over $S$ and such that $q^{-1}(U)$
is schematically dense.
\endit

\demo{Proof}  Indeed, if $X$ is proper over $S$ then in (\02.5.1) $X'$ is
proper over $S$, so the image of $X'$ in $\widebar X$ is closed.  Thus we
may assume that $\widebar X=X'$.\qed
\enddemo

\lemma{\02.7}  Let $X_1$ and $X_2$ be separated $S$\snug-schemes of finite
type, both of which contain open dense subschemes isomorphic to a scheme $U$.
Suppose that there exist quasi-dominations $\phi\:X_1\dashrightarrow X_2$
and $\psi\:X_2\dashrightarrow X_1$ compatible with the isomorphisms with $U$.
Then there exists a separated $S$\snug-scheme $X$ of finite type, and open
immersions $X_1\hookrightarrow X$ and $X_2\hookrightarrow X$, such that $U$ is
schematically dense in $X$, compatible with
$U\hookrightarrow X_i\hookrightarrow X$ for both $i$.
\endit

\demo{Proof}  The graph $\Gamma_1$ of $\phi$ is (by definition) a closed
subset of $X_1\times_S X_2$, isomorphic (via the first projection) to the
domain $W_1$ of $\phi$, and containing the image of $U$ in $X_1\times_S X_2$.
Since $U$ is dense in $X_1$, it follows that the set $\Gamma_1$ is the
closure of the image of $U$ in $X_1\times_S X_2$.  Likewise, the graph
$\Gamma_2\subseteq X_2\times_S X_1$ of $\psi$ is the closure of the image of $U$
in $X_2\times_S X_1$.  These two sets are therefore the images of each
other under the canonical isomorphism $X_1\times_S X_2\cong X_2\times_S X_1$.
Thus $\phi$ induces an isomorphism of $W_1$ with the domain $W_2$ of $\psi$,
and $\psi=\phi^{-1}$.

Let $X$ be the scheme obtained by glueing $X_1$ and $X_2$ along this
isomorphism.  Clearly the image of $U$ is a dense subscheme of $X$;
it remains only to show that $X$ is separated over $S$.  Let $\Delta$
denote the diagonal subset of $X\times_S X$.  The intersections of $\Delta$
with $X_1\times_S X_1$ and $X_2\times_S X_2$ are closed since $X_1$ and $X_2$
are separated over $S$.  The intersections of $\Delta$ with $X_1\times_S X_2$
and $X_2\times_S X_1$ are closed since $\phi$ and $\psi$ are quasi-dominations,
respectively.  Thus $\Delta$ is closed, so $X$ is separated over $S$,
as was to be shown.\qed
\enddemo

\prop{\02.8} \cite{N~1, 4.2}  Let
$$X_1 \hookleftarrow U \hookrightarrow X_2$$
be a diagram of separated $S$\snug-schemes of finite type, and assume
that $U$ is an open dense subscheme of both $X_i$.  Then there exists
a separated $S$\snug-scheme $X$ of finite type, and an open immersion
$U\hookrightarrow X$ with schematically dense image, such that $X$
properly quasi-dominates both $X_i$ compatible with the injections of $U$.
\endit

\demo{Proof}  By Theorem \02.4 there exists an ideal sheaf $\frak a_1$
on $X_1$, with $V(\frak a_1)$ disjoint from $U$, such that if
$\pi_1\:X_1'\to X_1$ denotes the blowing-up along $\frak a_1$, then $X_1'$
quasi-dominates $X_2$.  Likewise, there is an ideal sheaf $\frak a_2$ on $X_2$
such that the blowing-up $\pi_2\:X_2'\to X_2$ has the property that $X_2'$
quasi-dominates $X_1$.

We can go up another level as follows.  Let $p_1\:W_1\to X_2$ be the morphism
underlying the quasi-domination $X_1'\dashrightarrow X_2$, let $\frak a_1'$
be some extension of $p_1^{-1}\frak a_2\cdot\Cal O_{W_1}$ to $X_1'$,
let $\pi_1'\:X_1''\to X_1'$ denote the blowing-up of $X_1'$ along
$\frak a_1'$, and let $W_1'=(\pi_1')^{-1}(W_1)$.  Then, by
\cite{H, II Cor.~7.15}, $p_1$ extends uniquely to a morphism
$p_1'\:W_1'\to X_2'$.  We claim that this morphism determines a
quasi-domination $X_1''\dashrightarrow X_2'$; \ie, the graph $\Gamma_1'$
of $p_1'$ is closed in $X_1''\times_S X_2'$.  This graph is contained in
the closed subset $(\pi_1'\times_S\pi_2)^{-1}(\Gamma_1)$, where $\Gamma_1$
denotes the graph of $p_1$.  This proves the claim, since
$$\widebar{\Gamma_1'} \subseteq (\pi_1'\times_S\pi_2)^{-1}(\Gamma_1)
  \subseteq W_1'\times_S X_1',$$
and since $\Gamma_1'$ is a closed subset of $W_1'\times_S X_1'$.

Let $\frak a_2'$ and $\pi_2'\:X_2''\to X_2'$ be defined symmetrically;
then $X_2''\dashrightarrow X_1'$ is again a quasi-domination.

$$\xymatrix{
  X_1'' \ar[d]_{\pi_1'} \ar@{.>}[dr] & X_2'' \ar[d]^{\pi_2'} \ar@{.>}[dl] \\
  X_1' \ar[d]_{\pi_1} \ar@{.>}[dr] & X_2' \ar[d]^{\pi_2} \ar@{.>}[dl] \\
  X_1 & X_2 }$$

At this stage the process can be made to stop.  First, note that the graphs
of $p_1'$ and $p_2$ are the closures of the image of $U$ in $X_1''\times_S X_2'$
and $X_2'\times_S X_1$, respectively; therefore the image of $p_1'$ equals
the domain of $p_2$.  Hence
$$\split (p_1')^{-1}\frak a_2'\cdot\Cal O_{X_1''}
  &= (p_1')^{-1}\bigl(p_2^{-1}\frak a_1\cdot\Cal O_{X_2'}\bigr)
    \cdot\Cal O_{X_1''} \\
  &= (p_1')^{-1}p_2^{-1}\frak a_1\cdot(p_1')^{-1}\Cal O_{X_2'}
    \cdot\Cal O_{X_1''} \\
  &= (p_2\circ p_1')^{-1}\frak a_1\cdot\Cal O_{X_1''} \\
  &= (\pi_1\circ\pi_1')^{-1}\frak a_1\cdot\Cal O_{X_1''} \\
  &= (\pi_1')^{-1}\bigl(\pi_1^{-1}\frak a_1\cdot\Cal O_{X'}\bigr)
    \cdot\Cal O_{X''}.
\endsplit$$
But $\pi_1^{-1}\frak a_1\cdot\Cal O_{X'}$ is a line sheaf on $X_1'$, so we
may take $\frak a_1''$ to equal the pull-back of that line sheaf to $X_1''$.
Thus the blowing-up $X_1'''$ is isomorphic to $X_1''$.  Constructing
$\frak a_2''$ and $X_2'''$ the same way then leads to the situation of
Lemma \02.7.  Therefore there exists a scheme $X$ as in that lemma,
with $X_1''$ and $X_2''$ (and hence $U$) schematically dense in $X$.
Since $\pi_1\circ\pi_1'$ and $\pi_2\circ\pi_2'$ are proper, it follows
that $X$ properly quasi-dominates $X_1$ and $X_2$ (Lemma \02.2f).\qed
\enddemo

\beginsection{\03}{Some additional lemmas}

This section gives some additional lemmas that will be needed for the
main theorem in Section \04.

\lemma{\03.1}  Let $X$ be a scheme, separated and of finite type over a
noetherian scheme $S$.  Then there exists a cover $\{U_1,\dots,U_n\}$ of
$X$ by Zariski-dense open affines which are quasi-projective over $S$.
\endit

\demo{Proof}  Since $X$ is quasi-compact, it suffices to show that each
point $P\in X$ has an open neighborhood $U_P$ which is affine, Zariski-dense,
and quasi-projective over $S$.

First note that, since $X$ is of finite type over $S$, it is covered by
open affine subsets whose affine rings are of finite type over the affine
rings of open affine subsets of $S$.  In particular, these open subsets
are quasi-projective over $S$.  Each $P\in X$ has an open neighborhood of this
type; also, each irreducible component of $X$ contains an open subset
of this type disjoint from all other irreducible components
(take an open quasi-projective affine neighborhood of its generic point,
and then localize away from all other irreducible components of $X$).
(Since $X$ is noetherian, it has only finitely many irreducible components.)

Now fix $P\in X$.  Let $U_P^0$ be an open affine quasi-projective neighborhood
of $P$.  Let $U_P$ be the disjoint union of $U_P$ with open affine
quasi-projective subsets of all irreducible components of $X$ not
contained in $\widebar U_P^0$.  This is an open affine quasi-projective
Zariski-dense neighborhood of $P$.\qed
\enddemo

\lemma{\03.2}  Let
$$\CD X' @[hook]>>> \widebar X' \\
  @VVpV @VVqV \\
  X @[hook]>>> \widebar X\endCD$$
be a commutative diagram of schemes, such that $p$ and $q$ are separated
and of finite type, $q^{-1}(X)=X'$, $p$ is an isomorphism over an open
Zariski-dense subset $U$ of $X$, and the horizontal arrows are open immersions.
Let $F$ and $G$ be closed subsets of $X$ and $\widebar X'$, respectively,
such that $p^{-1}(F)\subseteq G$.

Then there exists an ideal $\frak a$ on $\widebar X$, with
$V(\frak a)\subseteq\widebar F\setminus F$, such that if
$\pi\:\widetilde X\to\widebar X$ denotes the blowing-up along $\frak a$,
if $\widetilde X'$ denotes the closure of $U$ in
$\widebar X'\times_{\widebar X}\widetilde X$, with morphisms as below,
$$\CD \widebar X' @<\pi'<< \widetilde X' \\
  @VVqV @VV\tilde qV \\
  \widebar X @<\pi<< \widetilde X\endCD\rlap{\quad,}\tag\03.2.1$$
and if $\widetilde F$ denotes the closure of $F$ in $\widetilde X$,
then $\tilde q^{-1}(\widetilde F)\subseteq(\pi')^{-1}(G)$.
\endit

\demo{Proof}  Apply Lemma \01.7 to the diagram
$$\CD \Gamma &\;\;\subseteq\;\;& (\widebar X'\setminus G)\times X
  @[hook]>>> (\widebar X'\setminus G)\times \widebar X \\
  && @VV\pr_2V @VV\pr_2V \\
  F &\;\;\subseteq\;\;& X @[hook]>>> \widebar X \endCD$$
with $\Gamma$ equal to the graph of the morphism
$p\restrictedto{\widebar X'\setminus G}$.  Let $\widetilde F$ denote as
usual the closure of $F$ in $\widetilde X$, and let $\widetilde\Gamma$
denote the closure of $U$ in $(\widebar X'\setminus G)\times\widetilde X$.
Then the lemma gives
$$\bigl((\widebar X'\setminus G)\times\widetilde F\bigr)
  \cap \widetilde\Gamma
  = \emptyset.$$
Since $\widetilde\Gamma\subseteq(\widebar X'\setminus G)\times\widetilde X$,
this can be shortened to
$(\widebar X'\times\widetilde F)\cap\widetilde\Gamma=\emptyset$.  Thus
$$\widetilde X'\cap(\widebar X'\times\widetilde F)
  \subseteq \widetilde X'\setminus\widetilde\Gamma\tag\03.2.2$$
(where $\widetilde X'$ is as in (\03.2.1)).

On the other hand, we have
$\widetilde\Gamma
  =\widetilde X'\cap\bigl((\widebar X'\setminus G)\times\widetilde X\bigr)$
(via the closed immersion
$\widebar X'\times_{\widebar X}\widetilde X
  \subseteq\widebar X'\times\widetilde X$); hence
$$\widetilde X'\setminus\Gamma' = \widetilde X'\cap(G\times\widetilde X)
  = (\pi')^{-1}(G)$$
(where again $\pi'$ is as in (\03.2.1)).  Also
$\tilde q^{-1}(\widetilde F)=\widetilde X'\cap(\widebar X'\times\widetilde F)$.
Thus, by (\03.2.2),
$$\tilde q^{-1}(\widetilde F) = \widetilde X'\cap(\widebar X'\times\widetilde F)
  \subseteq \widetilde X'\setminus\widetilde\Gamma = (\pi')^{-1}(G).\qed$$
\enddemo

\lemma{\03.3}  For $i=1,\dots,n$ let
$$\xymatrix@R=4mm{ & X_i \ar@{^{(}->}[]+<4mm,0mm>;[r] \ar[dd]^{q_i}
    & \widebar X_i \\ 
  U \ar@{^{(}->}[]+<4mm,2mm>;[ur] \ar@{^{(}->}[dr] \\ 
  & X }$$
be a diagram of separated $S$\snug-schemes of finite type,
with $U$ Zariski-dense in $X$, with $q_i$ proper and isomorphic over $U$,
with $X_i$ Zariski-dense in $\widebar X_i$,
and with $\widebar X_i$ proper over $S$.
Let $X^{*}$ be the closure of $U$ in the product of the $\widebar X_i$,
and let $p_i\:X^{*}\to\widebar X_i$ be the projection maps.
Finally, let $F_i$ be a closed subset of $X_i$ for each $i$,
such that $p_1^{-1}(F_1)\cap\dots\cap p_n^{-1}(F_n)=\emptyset$.

Then there exist ideal sheaves $\frak a_i$ on $\widebar X_i$ for each $i$,
with $V(\frak a_i)$ disjoint from $X_i$, such that after replacing each
$\widebar X_i$ with its blowing-up along $\frak a_i$, we obtain an
analogous situation in which
$$p_1^{-1}(\widebar F_1)\cap\dots\cap p_n^{-1}(\widebar F_n)=\emptyset.$$
\endit

\demo{Proof}  Let $\rho\:X^{**}\to X^{*}$ be the blowing-up of $X^{*}$
along the sum of the ideals associated to the sets $\widebar{p_i^{-1}(F_i)}$.
By Lemma \01.1, $\bigcap_{i=1}^n \widebar{\rho^{-1}(p_i^{-1}(F_i))}=\emptyset$.

Since the rational map $q_i\circ p_i\:X^{*}\to X$ does not depend on $i$,
and since $q_i$ is proper, the set $p_i^{-1}(X_i)$ does not depend on $i$,
and the closed set $\bigcap_{i=1}^n\widebar{p_i^{-1}(F_i)}$ does not meet
this open set.  Thus $\rho$ is an isomorphism over this open subset of $X^{*}$.
Therefore we may identify $p_i^{-1}(X_i)$ with an open subset of $X^{**}$.

Now apply Lemma \03.2 to the diagrams
$$\CD p_i^{-1}(X_i) @[hook]>>> X^{**} \\
  @VVV @VV p_i V \\
  X_i @[hook]>>> \widebar X_i\rlap{\quad,}\endCD$$
with $F=F_i$ and $G=\widebar{\rho^{-1}(p_i^{-1}(F_i))}$.
This provides us with the desired blowings-up.
Indeed, let $\widetilde X_1,\dots,\widetilde X_n$ denote the blowings-up
obtained from the lemma.  Let $X^{\#}$ be a proper $S$\snug-scheme that
dominates $X^{**}$ and the $\widetilde X_i$, with corresponding morphisms
$\phi\:X^{\#}\to X^{**}$ and $\psi_i\:X^{\#}\to\widetilde X_i$.
Finally, let $\widetilde F_i$ denote the closures
of the $F_i$ in $\widetilde X_i$.  We then have
$$\bigcap_{i=1}^n \psi_i^{-1}(\widetilde F_i)
  \subseteq \bigcap_{i=1}^n \phi^{-1}(\widebar{\rho^{-1}(p_i^{-1}(F_i))})
  = \phi^{-1}\left(\bigcap_{i=1}^n \widebar{\rho^{-1}(p_i^{-1}(F_i))}\right)
  = \emptyset.$$
Since the map from $X^{\#}$ to the new $X^{*}$ is surjective, the same holds
true on $X^{*}$.\qed
\enddemo

\beginsection{\04}{The main theorem}

It is now possible to prove Nagata's embedding theorem.

\thm{\04.1}  Let $X$ be a separated $S$\snug-scheme of finite type.
Then there exists an open immersion of $X$ into a proper $S$\snug-scheme,
with schematically dense image.
\endit

\demo{Proof}  Let $U_1,\dots,U_n$ be a cover of $X$ by open Zariski-dense
quasi-projective sets (Lemma \03.1).  For each $i$, Corollary \02.5 implies
that there exists a diagram
$$\xymatrix@R=4mm{ & X_i \ar@{^{(}->}[]+<4mm,0mm>;[r] \ar[dd]^{q_i}
    & \widebar X_i \\ 
  U_i \ar@{^{(}->}[]+<4mm,2mm>;[ur] \ar@{^{(}->}[dr] \\ 
  & X }$$
in which $q_i$ is a proper morphism, isomorphic over $U_i$, and in which
$\widebar X_i$ is proper over $S$.  Moreover, $U_i$ is schematically dense
in $\widebar X_i$.

Let $F_i=q_i^{-1}(X\setminus U_i)$, let $U=\bigcap_{i=1}^n U_i$,
let $X^{*}$ be the closure of $U$ in the product of the $\widebar X_i$,
and let $p_i\:X^{*}\to\widebar X_i$ be the restrictions of the projection
morphisms.  By Lemma \03.3, we may assume that the $\widebar X_i$ are
sufficiently blown up (without affecting $X_i$) so that
$$\bigcap_{i=1}^n p_i^{-1}(\widebar F_i) = \emptyset.$$

For each $i$ let $M_i$ be the $S$\snug-scheme obtained by glueing $X$
and $\widebar X_i\setminus\widebar F_i$ along $U_i$.  It is separated
because Lemmas \02.2f and \02.2b imply that $\widebar X_i\setminus\widebar F_i$
quasi-dominates $X$; hence the graph of the quasi-domination, which is
just the diagonal image of $U_i$, is closed in
$X\times_S(\widebar X_i\setminus\widebar F_i)$.  Also $X$ is schematically
dense in $M_i$ for all $i$, since $U_i$ is schematically dense
in $\widebar X_i$.

Let $M$ be the scheme obtained by applying Proposition \02.8 to the injections
$X\hookrightarrow M_i$ for all $i$.  Then $X$ is embedded as an open dense
subscheme of $M$.  It remains only to show that $M$ is proper over $S$.
Let $R$ be a valuation ring, let $K$ be its field of fractions, and let
$$\CD \Spec K @>>> M \\
  @VVV @VVV \\
  \Spec R @>>> S \endCD$$
be a commutative diagram.  Since $X^{*}$ is a proper $S$\snug-scheme, there is a
morphism $\Spec R\to X^{*}$ making a similar diagram commute.  Let $i$
be an index such that the image does not meet $p_i^{-1}(\widebar F_i)$.
Then there exists a morphism $\Spec R\to M_i$ making a similar diagram commute;
since $M$ properly quasi-dominates $M_i$, the same holds true for $M$.
This shows that $M$ is proper, by the valuative criterion of properness.\qed
\enddemo

\beginsection{\05}{Complements:  extensions of divisors and sheaves}

In number theory (at least), it is sometimes useful to be able to extend
an $S$\snug-scheme $X$ to a proper scheme over $S$ such that given
data extend to similar data on the larger scheme.  In many cases,
this follows by well-known techniques, and works for any given completion:
a sheaf of ideals extends by taking the scheme-theoretic closure of the
corresponding closed subscheme; a Weil divisor extends by taking the
closure; ditto for an effective Weil divisor; an arbitrary sheaf, sheaf
of $\Cal O_X$\snug-modules, or quasi-coherent sheaf extends by taking the
direct image; and a coherent sheaf extends by \cite{H, II Ex.~5.15}.

In other cases, one needs to be able to blow up the completion.  For
example, to extend an effective Cartier divisor $D$ to an effective Cartier
divisor on the completion, one can extend the ideal sheaf $\Cal O(-D)$
to an ideal sheaf on any given completion, blow up along the extended sheaf
to make the sheaf invertible, and then convert back to a Cartier divisor.

Extending arbitrary Cartier divisors and line sheaves is a bit more work,
but is still fairly easy, given the work that has already been done.
Extending vector sheaves involves further arguments, but still follows
the same plan as for Cartier divisors.  We will not treat line sheaves
separately, since they occur as a special case of vector sheaves, and
since the argument is parallel to that for Cartier divisors.

We start with a lemma used in extending Cartier divisors, then give a
series of lemmas leading up to a parallel result for vector sheaves,
and then prove the ultimate extension result for both objects simultaneously.

\lemma{\05.1}  Let $X$ be a separated $S$\snug-scheme of finite type,
and let $U$, $X_1$, and $X_2$ be open subschemes of $X$ with
$U\subseteq X_1\cap X_2$, $X_1\cup X_2=X$, and $U$ schematically dense in $X$.
Let $D_1$ and $D_2$ be Cartier divisors on $X_1$ and $X_2$, respectively,
that coincide on $U$.  Then there is a proper morphism $\pi\:X'\to X$
and a Cartier divisor $D'$ on $X'$, such that $\pi$ is an
isomorphism over $U$, such that $\pi^{-1}(U)$ is schematically dense in $X'$,
and such that the restriction of $D'$ to $\pi^{-1}(U)$ coincides with
the pull-backs of $D_1$ and $D_2$.
\endit

\demo{Proof}  Let $F=(X_1\cap X_2)\setminus U$.  After a preliminary
blowing-up, not affecting $U$, we may assume that the closure $\widebar F$
of $F$ in $X$ is the support of an effective Cartier divisor, also
called $\widebar F$ (see \cite{H, II~6.13a} and Lemma \01.3).
By Lemma \01.6, $U$ is still schematically dense in the new $X$.

We claim that, for sufficiently large $n$, $D_1-D_2+n\widebar F$
is an effective divisor on $X_1\cap X_2$.
By quasi-compactness it suffices to do this locally.  Let $\Spec A$ be
an open affine in $X_1\cap X_2$ such that $\widebar F$ is represented
by $f\in A$ and $D_1-D_2$ is represented by $g\in S^{-1}A$, where $S$
is the multiplicative system of elements of $A$ which are not zero divisors.
Since $D_1=D_2$ away from $F$, we actually have $g\in A_f^{*}$, and
the claim is then obvious.

It then follows that $\Cal O(-D_1+D_2-n\widebar F)$ is an ideal sheaf
on $X_1\cap X_2$; call it $\frak a$.  Extend it (by Lemma \01.2) to an
ideal sheaf on all of $X$, and let $\pi\:X'\to X$ be the blowing-up of $X$
along this ideal sheaf.  This blowing-up does not affect $X_1\cap X_2$,
and $\pi^{-1}(U)$ is schematically dense in $X'$.
The inverse image of the extended $\frak a$ is then equal to $\Cal O(-G)$
for some Cartier divisor $G$.  Therefore the Cartier divisors $D_1$ on $X_1$
and $D_2-n\widebar F+G$ on $X_2$ coincide on $X_1\cap X_2$, so they
combine to give a Cartier divisor $D'$ on $X'$, as was to be shown.\qed
\enddemo

This lemma will be applied in Proposition \05.6 and Theorem \05.7 (below);
the latter is the main result about extending Cartier divisors (and
vector sheaves) in the context of Nagata's embedding theorem.  First,
though, we prove a corresponding result for vector sheaves (in which
the rank-\snug$1$ case is very similar to the proof of Lemma \05.1,
but the higher-rank case introduces some additional difficulties).

\lemma{\05.2}  Let $B$ be a commutative ring, let $r\in\Bbb Z_{>0}$,
and let $N$ be a finitely generated submodule of $B^r$.  Assume that
the image of $\bigwedge^r N\to\bigwedge^rB^r\cong B$ is a principal ideal,
generated by the image of a decomposable element $n_1\wedge\dots\wedge n_r$,
and that this generator is a nonzerodivisor.  Then $N$ is a free $B$\snug-module
of rank $r$, with basis $n_1,\dots,n_r$.
\endit

\demo{Proof}  Let $\bold e_1,\dots,\bold e_r$ be the standard basis of $B^r$,
write $n_i=a_{i1}\bold e_1+\dots+a_{ir}\bold e_r$ for all $i$, and let
$A=(a_{ij})$ be the resulting $r\times r$ matrix.  The assumptions imply that
$\det A$ is a nonzerodivisor.  It then follows that $n_1,\dots,n_r$ are
linearly independent over $B$.  Indeed, if $b_1n_1+\dots+b_rn_r=0$ in $N$
with $b_1,\dots,b_r\in B$, and if $\vec b$ is the column vector composed
of $b_1,\dots,b_r$, then ${}^{\text t}A\vec b=\vec 0$.  Multiplying on the
left by the adjugate matrix of ${}^{\text t}A$ then gives
$(\det A)\vec b=\vec 0$, so $b_i=0$ for all $i$ since $\det A$ is a
nonzerodivisor.

To show that $n_1,\dots,n_r$ span $N$, let $n$ be an arbitrary element of $N$,
and for $i=1,\dots,r$ let $a_i\in B$ be the (unique) element
determined by the condition that the element
$$n_1\wedge\dots\wedge n \wedge\dots\wedge n_r - a_in_1\wedge\dots\wedge n_r$$
lies in the kernel of the map $\bigwedge^r N\to\bigwedge^rB^r$, where
in the first term the $n$ occurs in the $i^{\text{th}}$ position.
Let $n'=n-a_1n_1-\dots-a_rn_r$, and write $n'=b_1\bold e_1+\dots+b_r\bold e_r$
with $b_1,\dots,b_r\in B$.  Then all elements
$n_1\wedge\dots\wedge n' \wedge\dots\wedge n_r$ lie in the kernel of
$\bigwedge^r N\to\bigwedge^rB^r$, so if one replaces any row of the matrix
$A$ with $b_1,\dots,b_r$, then the determinant of the resulting matrix is zero.
(Here $A$ is the same matrix as in the preceding paragraph.)
Thus $b_1C_{i1}+\dots+b_rC_{ir}=0$ for all $i$, where $C_{ij}$ is the $(i,j)$
cofactor matrix of $A$.  This implies that
${}^{\text t}(\operatorname{adj}A)\vec b=\vec 0$; multiplying on the left
by ${}^{\text t}A$ gives $(\det A)\vec b=\vec 0$, and therefore $\vec b=\vec 0$.
Thus $n'=0$, and so $n$ lies in the span of $n_1,\dots,n_r$.\qed
\enddemo

\lemma{\05.3}  Let $X$ be a noetherian scheme, let $\Cal E$ be a vector
sheaf on $X$ of rank $r$, let $U$ be an open dense subscheme of $X$,
and let $\Cal F$ be a vector subsheaf of $\Cal E\restrictedto U$ (i.e.,
a locally free $\Cal O_U$\snug-submodule of $\Cal E\restrictedto U$),
also of rank $r$.  Assume that the induced map
$\bigwedge^r\Cal F\to\bigwedge^r\Cal E\restrictedto U$ is injective.
Then there exists a proper morphism $\pi\:X'\to X$, inducing an
isomorphism $\pi^{-1}(U)\to U$, and a vector subsheaf $\Cal F'$ of
$\pi^{*}\Cal E$ on $X'$ of rank $r$, such that
$\Cal F'\restrictedto{\pi^{-1}(U)}=\pi^{*}\Cal F$, and such that $\pi^{-1}(U)$
is schematically dense in $X'$.
\endit

\demo{Proof}  The result is trivial if $r=0$, so we assume that $r>0$.

Following \cite{EGA, I~9.5.2}, we let
$\Cal G
 =\ker\bigl(\Cal E\to i_{*}\bigl(\Cal E\restrictedto U/\Cal F\bigr)\bigr)$,
where $i\:U\to X$ is the inclusion map.
This is a coherent sheaf by \cite{H, II~5.8c}, \cite{H, II~5.7}, and the fact
that it is a subsheaf of the coherent sheaf $\Cal E$.
The image of $\bigwedge^r\Cal G$ in $\bigwedge^r\Cal E$ is a coherent subsheaf
of a line sheaf, so it defines an ideal sheaf $\frak a$ on $X$.
Let $\pi\:X'\to X$ be the blowing-up of $X$ along $\frak a$.
Since $\Cal G\restrictedto U=\Cal F$ and since
$\bigwedge^r\Cal F\to\bigwedge^r\Cal E$ is injective, it follows
that $\frak a\restrictedto U$ is a line sheaf, and therefore $\pi$ is
an isomorphism over $U$.  Let $\Cal F'$ be the image of $\pi^{*}\Cal G$
in $\pi^{*}\Cal E$.  Since $\Cal G\restrictedto U=\Cal F$, we have
$\Cal F'\restrictedto{\pi^{-1}(U)}=\pi^{*}\Cal F$.

It remains only to show that $\Cal F'$ is a vector sheaf of rank $r$.
Let $\Spec A$ be an open affine in $X$ over which $\Cal E$ is free,
fix an isomorphism $\Gamma(\Spec A,\Cal E)\cong A^r$, and let $\Spec B$
be an open affine in $\pi^{-1}(\Spec A)$.  Let $M$ be the submodule of $A^r$
corresponding to $\Cal G\restrictedto{\Spec A}$, and let $N$ be the image
of $M\otimes_A B$ in $B^r$, so that $\Cal F'\restrictedto{\Spec B}=\widetilde N$
(as subsheaves of $\pi^{*}\Cal E\restrictedto{\Spec B}=\widetilde{B^r}$).
The ideal $\Gamma(\Spec A,\frak a)$ in $A$ corresponds to the image of
$\bigwedge^r M \to\bigwedge^r A^r\cong A$; let $\frak b$ denote the image
of this ideal in $B$.  By \cite{H, II~7.13a}, $\widetilde{\frak b}$ is a
line sheaf, so after localizing $B$ further we may assume that $\frak b$
is a principal ideal, generated by a nonzerodivisor.

We note that the submodule $\frak b_0$ of $\bigwedge^rB^r$ corresponding
to $\frak b$ is generated by elements of the form
$(m_1\wedge\dots\wedge m_r)\otimes_A 1$ with $m_1,\dots,m_r\in M$.
On the other hand, the image of $\bigwedge^r N$ in $\bigwedge^r B^r$
is generated by elements of the form
$(m_1\otimes1)\wedge\dots\wedge(m_r\otimes1)$ for $m_1,\dots,m_r\in M$,
and these elements correspond to the above generators for $\frak b_0$
under the isomorphism
$\bigl(\bigwedge^r A^r\bigr)\otimes_A B \cong \bigwedge^r(A^r\otimes_A B)$.
This shows that the image of $\bigwedge^r N\to\bigwedge^r B^r\cong B$ also
equals $\frak b$, and therefore is principal, generated by a nonzerodivisor.

We next claim that one such (single) generator can be written as a decomposable
element $(m_1\wedge\dots\wedge m_r)\otimes1$, assuming that $B$ has been
chosen properly.  To show this, we note that the ideal
$\frak a^{*}:=\Gamma(\Spec A,\frak a)$ is generated by elements of $A$
of the form $\rho(m_1\wedge\dots\wedge m_r)$,
where $\rho\:\bigwedge^rA^r\overset\sim\to\to A$ is the usual isomorphism
and $m_1,\dots,m_r\in M$.  For any given choice of $m_1,\dots,m_r$, let
$x=\rho(m_1\wedge\dots\wedge m_r)$, and let
$B=\bigl(\bigoplus_{d\ge 0}(\frak a^{*})^d\bigr)_{(x)}$, with $x$ of
degree $1$.  Then $\Spec B\subseteq\pi^{-1}(\Spec A)$,
and $\frak b=\frak a^{*}B$ is a principal ideal, generated by $x$
(in degree $0$).  We may then apply Lemma \05.2 to conclude that $\Cal F'$
is free of rank $r$ over $\Spec B$.  As $m_1,\dots,m_r$ vary, the
corresponding sets $\Spec B$ cover $\pi^{-1}(\Spec A)$, and thus $\Cal F'$
is locally free on $X'$ of rank $r$.\qed
\enddemo

\remk{\05.4}  This proof also shows that
$\bigwedge^r\Cal F'\to\bigwedge^r\pi^{*}\Cal E$ is injective (for the given
choice of $\Cal F'$).
\endit

\lemma{\05.5}  Let $X$ be a separated $S$\snug-scheme of finite type,
and let $U$, $X_1$, and $X_2$ be open subschemes of $X$ with
$U\subseteq X_1\cap X_2$, $X_1\cup X_2=X$, and $U$ schematically dense in $X$.
Let $\Cal E_1$ and $\Cal E_2$ be vector sheaves on $X_1$ and $X_2$,
respectively, such that $\Cal E_1\restrictedto U\cong\Cal E_2\restrictedto U$.
Then there is a proper morphism $\pi\:X'\to X$ and a vector sheaf $\Cal E'$
on $X'$, such that $\pi$ is an isomorphism over $U$, such that $\pi^{-1}(U)$
is schematically dense in $X'$, and such that
$\Cal E'\restrictedto{\pi^{-1}(U)}
 \cong\pi^{*}\bigl(\Cal E_1\restrictedto U\bigr)$.
\endit

\demo{Proof}  By treating each connected component of $U$ separately, we may
assume that the $\Cal E_i$ have constant rank $r$.

After a preliminary blowing-up, not affecting $U$, we may assume that
$X\setminus U$ is the support of an effective Cartier divisor $D$.
Let $\phi\:\Cal E_1\restrictedto U\to\Cal E_2\restrictedto U$ be the
given isomorphism.  For any open affine $\Spec A$ in $X_1\cap X_2$ on which
$D$ is principal, given say by $(f)$ with $f\in A$, and on which $\Cal E_1$
and $\Cal E_2$ are trivial, the morphism $\phi$ is given by an invertible
$r\times r$ matrix $M$ with entries in $A_f$.  Since $D$ is Cartier, $f$
is a nonzerodivisor, the map $A\to A_f$ is injective, and therefore for
all sufficiently large integers $n$ the matrix $f^nM$ extends uniquely
to a matrix with entries in $A$.  By quasi-compactness, there is one
integer $n$ which works for all such open affines $\Spec A$.  After replacing
$\Cal E_2$ with $\Cal E_2\otimes\Cal O(nD)$, we may assume that the isomorphism
$\phi$ extends (uniquely) to an injective morphism
$\Cal E_1\restrictedto{X_1\cap X_2}\to\Cal E_2\restrictedto{X_1\cap X_2}$
of $\Cal O_{X_1\cap X_2}$\snug-modules.

Now let $\Spec A$, $f$, and $M$ be as above.  Since the entries of $M^{-1}$
also lie in $A_f$, there is an integer $m$ such that the entries of $f^mM^{-1}$
all lie in $A$, and therefore $\det M\mid f^{rm}$ in $A$.  In particular,
$\det M$ is a nonzerodivisor, so letting $\Spec A$ vary over an open cover
of $X_1\cap X_2$, we see that $\bigwedge^r\phi\restrictedto{X_1\cap X_2}$
is injective.  By Lemma \05.3, there is then a proper morphism
$\pi_0\:X_2'\to X_2$, isomorphic over $X_1\cap X_2$, and a vector subsheaf
$\Cal E_1'$ of $\pi_0^{*}\Cal E_2$ extending the pull-back of $\Cal E_1$ on
$X_1\cap X_2$.  We then glue the $X$\snug-schemes $X_2'$ and $X_1$ over
$X_1\cap X_2$ to obtain a proper morphism $\pi\:X'\to X$, and glue the
vector sheaves $\Cal E_1$ on $X_1$ and $\Cal E_1'$ on $X_2'$ to obtain
the desired vector sheaf $\Cal E'$ on $X'$.  Since $U$ is schematically
dense in $X_1\cap X_2$ and $\pi^{-1}(X_1\cap X_2)$ is schematically dense
in $X'$, it follows that $\pi^{-1}(U)$ is schematically dense in $X'$.\qed
\enddemo

\prop{\05.6}  Let
$$X_1 \hookleftarrow U \hookrightarrow X_2$$
be as in Proposition \02.8, and let $D$ be a Cartier divisor (resp\.
let $\Cal E$ be a vector sheaf) on $U$.  Assume that $D$ (resp\. $\Cal E$)
extends (separately) to Cartier divisors (resp\. vector sheaves) on $X_1$
and on $X_2$.  Then one may choose $X$ in Proposition \02.8 so that $D$
(resp\. $\Cal E$) extends to a Cartier divisor (resp\. vector sheaf) on $X$.
\endit

\demo{Proof}  For $i=1,2$ let $D_i$ be a Cartier divisor on $X_i$ that
extends $D$ (resp\. let $\Cal E_i$ be a vector sheaf on $X_i$ that extends
$\Cal E$).

By Proposition \02.8 there exists a separated $S$\snug-scheme $X_0$
of finite type, together with an inclusion $U\hookrightarrow X_0$,
such that $X_0$ properly quasi-dominates $X_1$ and $X_2$, compatible with
the inclusions of $U$ into $X_0$, $X_1$, and $X_2$.  Replace $X_1$ and $X_2$
with the domains of these quasi-dominations, and $D_1$ and $D_2$
(resp\. $\Cal E_1$ and $\Cal E_2$) with their pull-backs, so that $X_1$
and $X_2$ are open dense subschemes of $X_0$.  After shrinking $X_0$,
we may assume that $X_0=X_1\cup X_2$.

We then conclude by applying Lemma \05.1 (resp\. Lemma \05.5).\qed
\enddemo

\thm{\05.7}  Let $X$ be a separated $S$\snug-scheme of finite type.
Then, given any Cartier divisor $D$ (resp\. vector sheaf $\Cal E$)
on $X$, the immersion $X\hookrightarrow\widebar X$ of Theorem \04.1
can be chosen so that $D$ (resp\. $\Cal E$) extends to a Cartier divisor
(resp\. vector sheaf) on $\widebar X$.
\endit

\demo{Proof}  This will proceed by making just a few changes to the
proof of Theorem \04.1.  We use the notation of that proof.

The open cover $U_1,\dots,U_n$ may be chosen so that $D\restrictedto{U_i}$
is principal for each $i$, say $D\restrictedto{U_i}=(f_i)$ (resp\.
so that $\Cal E$ is trivial on $U_i$).  Then $D$ extends to $M_i$
by requiring that $D\restrictedto{\widebar X_i\setminus\widebar F_i}=(f_i)$
(since $U_i$ is schematically dense in $\widebar X_i\setminus\widebar F_i$)
(resp\. $\Cal E$ extends as a trivial vector sheaf).  The rest of the proof
continues as before, but with Proposition \05.6 replacing Proposition \02.8.\qed
\enddemo

\remk{\05.8}  When $X$ is regular, the ability to extend a vector sheaf
$\Cal E$ has already been noted \cite{C-G, p.~23}.
\endit

\thm{\05.9}  Let $X$ be a separated $S$\snug-scheme of finite type.
Then, given any finite collections $D_1,\dots,D_r$ of Cartier divisors
and $\Cal E_1,\dots,\Cal E_s$ of vector sheaves on $X$, the immersion
$X\hookrightarrow\widebar X$ of Theorem \04.1 can be chosen so that
$D_1,\dots,D_r$ and $\Cal E_1,\dots,\Cal E_s$ extend to $\widebar X$
as Cartier divisors and vector sheaves, respectively.
\endit

\demo{Proof}  First note that if $r=s=0$ then this is just Theorem \04.1,
so we may assume that $r>0$ or $s>0$.

Let $i_j\:X\hookrightarrow\widebar X_j$, $j=1,\dots,r+s$, be completions
for which $D_1,\dots,D_r,\Cal E_1,\dots,\Cal E_s$ extend, respectively.
Then we may let $\widebar X$ be the scheme-theoretic image of the map
$$(i_1,\dots,i_{r+s})\:X
  \hookrightarrow \widebar X_1\times_S\dotsm\times_S\widebar X_{r+s}$$
and pull the extended data back via the respective projection morphisms.\qed
\enddemo

\remk{\05.10}  In Theorem \05.9, $X$ is dense in the resulting completion
$\widebar X$, so if any of the $\Cal E_i$ have constant rank, then so do
their extensions to $\widebar X$.
\endit

Finally, we note that the extension of Chow's lemma given here leads to
an extension of Theorem \05.9 in the (quasi-)projective setting.

\cor{\05.10}  In Theorem \05.9, if $X$ is quasi-projective over $S$,
then $\widebar X$ can be taken to be projective over $S$.
\endit

\demo{Proof}  Apply Corollary \02.6 to a completion
$X\hookrightarrow\widetilde X$ as in Theorem \05.9, and pull back the
Cartier divisors and vector sheaves to the resulting projective scheme
$\widebar X$.\qed
\enddemo

\comment
\beginsection{Exercises}

1.  Verify the isomorphism $(T_i)_{(b_j/a_i)}\cong S'_{(a_ib_j)}$ in the
proof of Proposition \01.5.

2.  Prove that, in Lemma \05.5, the morphism $\pi$ can be taken to be the
blowing-up of a sheaf of ideals $\frak a$ on $X$ with $V(\frak a)$ disjoint
from $U$.  ||| Can I go further with this?

\endcomment

\enddocument